\documentclass[12pt]{amsart}
\usepackage{amsmath,amstext,amscd}
\usepackage{amssymb}
\usepackage{hyperref}

\newtheorem{defin}{Definition}[section]
\newtheorem{teor}[defin]{Theorem}
\newtheorem{corol}[defin]{Corollary}
\newtheorem{osserv}[defin]{Remark}
\newtheorem{lemma}[defin]{Lemma}
\newtheorem{propos}[defin]{Proposition}
\newcommand {\n}{\noindent}

\begin{document}


\title{Finite solvable groups whose Quillen
       complex is Cohen-Macaulay}

\author{Francesco Matucci}
\address{Department of Mathematics, Cornell University, \\
310 Malott Hall, Ithaca, NY 14853, USA}
\email{matucci@math.cornell.edu}



\begin{abstract}
We will prove that the $p$-Quillen complex of a
finite solvable group with cyclic derived group is Cohen-Macaulay, if
$p$ is an odd prime. If $p = 2$ we prove a similar conclusion, but there is
a discussion to be made.
\end{abstract}

\maketitle

\markboth{Finite solvable
groups whose Quillen complex is Cohen-Macaulay}{Francesco Matucci}


\section{Introduction}

For a finite group $G$ and a prime number $p$, following
\cite{qu}, we denote by
$\mathcal{S}_p(G)$ the partially ordered set (in short: poset)
of all non-trivial $p$-subgroups of $G$ and by $\mathcal{A}_p(G)$ the poset
of all non-trivial $p$-elementary abelian subgroups of $G$,
both ordered by inclusion. We denote by
$\Delta(\mathcal{S}_p(G))$ and $\Delta(\mathcal{A}_p(G))$ the associated
chain complexes (or rather their geometric realizations); these are respectively called,
the $p$-\emph{Brown complex} and the
$p$-\emph{Quillen complex} of $G$, and,
as Quillen pointed out in \cite{qu},
they are homotopically equivalent.

One of the motivations for studying the Quillen complex
is a famous conjecture of Quillen himself which says that
$\Delta(\mathcal{A}_p(G))$ is contractible if and only if $G$ possesses
a non-trivial normal $p$-subgroup. This conjecture has been proved
by Quillen himself in various cases (in particular for solvable
groups); the most general known result on this conjecture is due to
M. Aschbacher and S. Smith (see \cite{assm}).

Quillen's conjecture is just one aspect of the more general question of
describing the homotopy type of the Quillen complex of a group G and its
connections with the algebraic structure of $G$. In the same seminal paper
(\cite{qu}), Quillen proved the following result.

\begin{teor}[Quillen]\label{tq}
Let $p$ be a prime number and $G=NP$ be a solvable group such that
$N \trianglelefteq G$ is a $p'$-group and $P$ is a $p$-elementary
abelian group. Then $\Delta(\mathcal{A}_p(G))$ is Cohen-Macaulay.
\end{teor}

Before continuing, let us explain some of the terminology
we adopt, as it might not be fully standard. 
Let $\Gamma$ be a finite simplicial complex, and $r\ge 0$ an integer.
We say that $\Gamma$ is {\sl weakly $r$--spherical},
if $\widetilde{H}_q(\Gamma) = 0$,
for all $q \ne r$; while we say that $\Gamma$ is $r$--{\sl spherical}
if it has the same homotopy type of a wedge of $r$-spheres or
it is contractible. We say that $\Gamma$ is {\sl spherical} (respectively,
{\sl weakly spherical}) if there exists an integer $r \ge 0$ such that
$\Gamma$ is $r$-spherical (respectively, weakly $r$-spherical).

Let $\sigma$ be a simplex of $\Gamma$. The {\sl link}
of $\sigma$ in $\Gamma$ is the subcomplex
$L_{\Gamma}(\sigma) = \{\tau \in \Gamma \mid \tau \cup \sigma \in \Gamma,
\tau \cap \sigma = \varnothing \}$. Let $d$ be the dimension of $\Gamma$:
finally we say
that $\Gamma$ is \textbf{Cohen-Macaulay}
provided $\Gamma$ is $d$-spherical and  the link of
each $r$-simplex of $\Gamma$ is $(d-r-1)$-spherical.

From Theorem \ref{tq} it easily follows that the $p$-Quillen complex of a
solvable group whose Sylow $p$-subgroups are abelian is spherical
(and indeed Cohen-Macaulay).
In this paper our aim is to give a description
of the $p$-Quillen complex in the more general case of solvable
groups whose Sylow $p$-subgroups have a cyclic derived group. 

\begin{teor}\label{main}
Let $G$ be a finite solvable group and $p$ an odd prime number dividing the order of $G$.
Let $P \in \text{\emph{Syl}}_p(G)$ and suppose that $P'$ is cyclic. Then $\Delta(\mathcal{A}_p(G))$
is Cohen-Macaulay.
\end{teor}

We point out that recently Fumagalli \cite{fum} has shown that, for $p$ odd,
the $p$-Quillen complex of a solvable group always has the homotopy type
of a wedge of spheres (of possibly different dimensions). Also, we observe
that, to the best of our knowledge, it is not known whether or not Quillen's result \ref{tq}
holds without the assumption of the solvability of $N$. 

For the case $p=2$,
we may have a different behavior than the one described in the previous Theorem.
We describe the relevant structure of
Sylow 2-subgroups and we show that it is possible to build
examples which are not even weakly spherical.

\begin{teor}\label{example}
Let $G$ be a finite solvable group such that $2$ divides its order.
Let $P \in \text{\emph{Syl}}_2(G)$ and suppose that $P'$ is cyclic. Then 
\medskip \\
(i) $\Omega_1(P)=TD$, $T,D \trianglelefteq \Omega_1(P)$, $\vert T \cap D \vert \le 2$,
$D = Z(D)E$, with $E=1$ or $E$ extra-special, and  either $T=1$ or
$T$ is dihedral or $T$ is semi-dihedral.
\medskip \\
(ii) Suppose that $\Omega_1(P)$ is the central product of $T$ and $D$. 
If $T = 1$ or $T$ is dihedral then $\Delta(\mathcal{A}_p(G))$ is Cohen-Macaulay.
If $T$ is semi-dihedral then there exist two distinct reduced homology groups of
$\Delta(\mathcal{A}_p(G))$ which are non-trivial.
\end{teor}

Let us now fix some notations. If $r$ is an element of a poset $P$ we denote
with $P_{>r}$ the poset $\{q \in P \mid q > r \}$ and say the
$P_{>r}$ is an upper interval. Analogously we define $P_{<r}$,
$P_{\ge r}$ and $P_{\le r}$.
Moreover if $r < s$ are elements of a poset $P$, we denote with
$(r,s)$ the poset $\{q \in P \mid r < q < s \}$.
If $\Gamma_1$ and $\Gamma_2$ are two complexes, we define \emph{wedge} of
$\Gamma_1$ and $\Gamma_2$ (denoted by $\Gamma_1 \vee_{x=y}
\Gamma_2$) to be the complex whose vertex set is given by the disjoint union of the vertices of $\Gamma_1$ and those
of $\Gamma_2$, by identifying a vertex $x$ of $\Gamma_1$ with a vertex $y$ of $\Gamma_2$ and whose set of simplices is
given by the union of the simplices of $\Gamma_1$ and those of $\Gamma_2$. We warn the reader that, as in
\cite{pw}, in our wedge decompositions of $\Delta(\mathcal{A}_{p}(G))$
the wedge of spaces is not formed using just a single wedge point, instead we
always have to specify where it is wedged in for each space.
We define the \emph{join} of the
simplicial complexes $\Gamma_1$ and $\Gamma_2$ (denoted by $\Gamma_1 \ast \Gamma_2$) 
to be the simplicial complex whose vertex set is
given by the disjoint union of the vertices of $\Gamma_1$ and those of $\Gamma_2$ and whose set of simplices is
given by the union of the simplices of $\Gamma_1$ and those of $\Gamma_2$ and the set with elements given by the
disjoint unions $\sigma \sqcup \tau$ of simplices $ \sigma \in \Gamma_1$ and $\tau \in \Gamma_2$. 
For every complex $\Gamma$ we define the join $\Gamma \ast \emptyset$ to be
equal to the complex $\Gamma$ itself.
We term $p$-\emph{torus} every $p$-elementary abelian group,
and we call the
\emph{rank} of a $p$-group $P$ (denoted by $rk(P)$), the maximal
dimension of a $p$-torus of $P$ as a $\mathbb{Z}_p$-vector
space. If $P$ is a $p$-group we denote by $\Omega_1(P)$ the subgroup
of $P$ generated by the elements of $P$ of order $p$.
For the rest we follow the notation of \cite{go}.


\section{Topological tools and first reductions}

We start by collecting a few topological results which we will use later.
Since we work with simplicial complexes in finite groups, we
will assume that all of the complexes will have finitely many simplices.

\begin{lemma}[Gluing Lemma]\label{gl}
Let $X$ be a simplicial complex, $Y_1, \ldots, Y_k,\ Z$ subcomplexes of
$X$ such that $X=Y_1 \cup \ldots \cup Y_k$,
$\dim X=\dim Y_i=r$ and $(r-1) \le \dim Z \le r$.
Suppose that $Z \subseteq Y_i$, for all $i=1, \ldots, k$, and that
$Y_i \cap Y_j = Z$, for all $i \ne j$. 
\medskip \\
(i) If $Z$ is weakly $(\dim Z)$-spherical and
$Y_i$ is weakly $r$-spherical for all $i$, then $X$ is weakly $r$-spherical.
\medskip \\
(ii) If $Z$ is $(\dim Z)$-spherical and
$Y_i$ is $r$-spherical for all $i$, then $X$ is $r$-spherical.
\end{lemma}

\noindent \emph{Proof.} (i) By induction on $k$. If $k=1$ there is nothing
to be proven. Let $k > 1$ and $Y=Y_1 \cup \ldots \cup Y_{k-1}$. By
induction on $Y$ we have that $Y$ is weakly $r$-spherical.
Observing that $Y \cap Y_k=
(Y_1 \cap Y_k) \cup \ldots \cup (Y_{k-1} \cap Y_k)=Z$,
by applying the Mayer-Vietoris sequences we get
$$
\ldots \to \widetilde{H}_{q}(Y \cap Y_{k}) \to \widetilde{H}_{q}(Y)
\oplus \widetilde{H}_{q}(Y_{k}) \to \widetilde{H}_{q}(X) \to
\widetilde{H}_{q-1}(Y \cap Y_{k}) \to \ldots $$

Now if $q \le (r-1)$ we have that
$0 = 0 \oplus 0 =
\widetilde{H}_{q}(Y) \oplus \widetilde{H}_{q}(Y_{k}) \to
\widetilde{H}_{q}(X) \to
\widetilde{H}_{q-1}(Y \cap Y_{k}) = \widetilde{H}_{q-1}(Z)=0$ and
so $\widetilde{H}_{q}(X) = 0$ by the exactness of the sequence. 

(ii) If $r=0$, it is obvious. If $r=1$, then
the $1$-dimensional complex $X$ is a graph and therefore it
is known to be contractible or $1$-spherical, for it is a union of segments
and circuits and hence it is connected. So we can suppose that $r \ge 2$.
By part (i) we have that $X$ is weakly $r$-spherical;
so, by the Hurewicz-Whitehead Theorem (see \cite{sp}),
it remains to be proven that $\pi_1(X)=1$. This can be achieved by applying the Van Kampen's
Theorem for simplicial complexes (see Theorem 11.60 in \cite{rotman}):
the Theorem states that $\pi_1(X) \cong (\pi_1(Y_1) \ast \ldots \ast \pi_1(Y_k)/\sim)$ for a suitable equivalence relation $\sim$.
Since $\pi_1(Y_i)=1$, for all $i$, this implies that $\pi_1(X) \cong (1 \ast \ldots \ast 1 /\sim) \cong 1$. $\square$

Let us now turn to groups.
A key tool in our proofs is the following result due to Pulkus and Welker,
which has been proven in \cite{pw}.

\begin{teor}[Pulkus-Welker]\label{PW}
Let $G$ be a finite group with a solvable normal $p'$-subgroup $N$. For
$A \le G$ set $\overline{A}=AN/N$. Then $\Delta(\mathcal{A}_{p}(G))$ is
homotopically equivalent to the wedge
\[
\Delta(\mathcal{A}_{p}(\overline{G})) \vee
\underset{\overline{A}
\in \mathcal{A}_{p}(\overline{G})}
\bigvee
\Delta(\mathcal{A}_{p}(NA)) \ast
\Delta(\mathcal{A}_{p}(\overline{G})_{>\overline{A}})
\]
where for each $\overline{A} \in \mathcal{A}_{p}(\overline{G})$ an arbitrary
point $c_{\overline{A}} \in \Delta(\mathcal{A}_{p}(NA))$ is identified with
$\overline{A} \in \Delta(\mathcal{A}_{p}(\overline{G}))$.
\end{teor}

We recall that, by using standard topological methods, one
proves that the wedge of many
$r$-spherical complexes is again an $r$-spherical
complex, and that the join of an $r$-spherical complex with an $s$-spherical
complex is an $(r+s+1)$-spherical complex.

The proof of our main result proceeds by first dealing with a split case
as in Theorem \ref{tq}. Thus,
fix a prime $p$ and consider a semidirect product  $G = N \rtimes P$,
where $P \ne 1$ is a $p$-group and $N \ne 1$ is a solvable $p'$-group.

\begin{lemma} \label{23}
Let $p$ be a prime number and $G=NP$ be a solvable group such that
$N \trianglelefteq G$ is a $p'$-group and $P$ is a $p$-group.
Suppose that, for any $X \in \mathcal{A}_p(P)$, the complex
$\Delta(\mathcal{A}_p(P)_{>X})$ is $(rk(P)-rk(X)-1)$-spherical. Then
the complex $\Delta(\mathcal{A}_p(G))$ is $(rk(P)-1)$-spherical.
\end{lemma}

\noindent \emph{Proof.} We apply the Pulkus-Welker formula to the complex
$\Delta(\mathcal{A}_p(G))$ and see that:
\[
\Delta(\mathcal{A}_{p}(G)) \simeq
\underset{X \in \mathcal{A}_{p}(P)}
\bigvee
\Delta(\mathcal{A}_{p}(NX)) \ast
\Delta(\mathcal{A}_{p}(P)_{>X}).
\]
Let $X \in \mathcal{A}_p(P)$. Then, 
by Quillen's Theorem \ref{tq}, we have that $\Delta(\mathcal{A}_p(NX))$ is $(rk(X)-1)$-spherical.
By hypothesis, we have that $\Delta(\mathcal{A}_{p}(P)_{>X})$ is $(rk(P)-rk(X)-1)$-spherical and hence that
$\Delta(\mathcal{A}_{p}(NX)) \ast \Delta(\mathcal{A}_{p}(P)_{>X})$ is spherical of rank
$(rk(X)-1)+(rk(P)-rk(X)-1)+1=rk(P)-1$. As we have observed before, the wedge of complexes that
are $(rk(P)-1)$-spherical is again $(rk(P)-1)$-spherical. $\square$

This reduces us to study the behavior of all the upper intervals
$\mathcal{A}_{p}(P)_{>X}$, for a non trivial $p$-group $P$
(in our case $P$ will have a cyclic derived group),
 and $X \in \mathcal{A}_{p}(P)$.

\noindent \textbf{Notation.} In the following we
will often drop the notation $\Delta(P)$ for a
poset $P$, and we will refer to $P$ both in the case of a poset and
in the case of a simplicial complex.


\section{$p$-groups with $P=\Omega_1(P)$ and $P'$ cyclic}

Let $P$ be a $p$-group, and $X\in \mathcal{A}_{p}(P)$.
Observe that $$\mathcal{A}_{p}(P)_{>X} =
\mathcal{A}_{p}(\Omega_1(P))_{>X},$$
so we may well assume that $\Omega_1(P)=P$.

\n We need to describe the $p$-groups $P$ such that $P=\Omega_1(P)$
and $P'$ is cyclic. The following results are certainly known, but
we include proofs for completeness. If $p>2$ such groups are
 ``essentially''  extra-special groups.
\begin{teor}\label{31}
Let $p$ be an odd prime number and $P$ a $p$-group such that $P'$ is cyclic
and $\Omega_1(P)=P$. If $P'\neq 1$, then $\vert P' \vert = p$, and $P$
is the direct product of an elementary abelian group and an extra-special
group of exponent $p$.
\end{teor}
\noindent \emph{Proof.}
Since $P'$ is cyclic, $P$ is regular, by Theorem 4.3.13 in \cite{su2}.
Now we use Theorem 4.3.14 in \cite{su2}
to obtain that $\exp(\Omega_1(P))=p$. But
$\Omega_1(P)=P$ and so $\vert P' \vert \le p$. The last assertion follows
easily. \ $\square$

For $p=2$ the situation is more complicated. For instance, dihedral and
semidihedral $2$-groups have cyclic derived subgroups and are generated by
involutions.
\begin{teor}\label{32}
Let $P$ be a $2$-group such that $\Omega_1(P)=P$, $P'$ is cyclic and
$\vert P' \vert > 2$. Then $P=TD$, with $T,D \trianglelefteq P$,
$\vert T \cap D \vert = 2$, $D = Z(D)E$, $E = 1$ or $E$ is extra-special
and $T$ is dihedral or $T$ is semi-dihedral. Moreover $\Omega_1(Z(P)) \le D$.
\end{teor}
\noindent \emph{Proof.} Since $P$ and all of its images are generated by
involutions we have $P' = \Phi (P) = P^2$, where
$P^2 = \langle x^2 \mid x \in P \rangle$. As $P'$ is cyclic, we have in
particular $P'= \langle x^2\rangle$ for some $x\in P$. Since $|P'|>2$, 
$| x | = 2^n$, for some $n\geq 3$.

We write: $X = \langle x\rangle$, $Z = \Omega_1(X) = \langle
x^{2^{n-1}}\rangle$, and $C = C_P(P')$; note that all these subgroups are
normal in $P$, $Z$ is central, and $Z<P'$.

Now, $C \geq C_P(X) \geq \Phi (P)$, hence $P/C_P(X)$ embeds in
$\Omega_1(Aut (X))$, which is an elementary abelian group of order $4$.
If $\vert P/C_P(X) \vert = 4$, there exists a $g\in P$ such that
$x^g = x^{1+2^{n-1}}$; but then $(x^2)^g = x^2$, and so $g\in C\setminus
C_P(X)$. Thus, in any case, $\vert P/C\vert \leq 2$.

Now let $r,s \in C$; then $[r,s,s]= 1$, and so
$[r,s]^2 = [r^2,s]=1$, since $r^2\in P'$. Therefore $C'$ is
generated by involutions; as $C'\leq P'$ is cyclic, we conclude that
$\vert C' \vert \le 2$ (and $C'\leq Z)$. By assumption, $\vert
P'\vert \geq 4$, thus we have $C < P$ (and $\vert P/C\vert = 2$).

From $C<P$ and $P= \Omega_1(P)$, we conclude that there exists an
involution $z\in P\setminus C$. Let $T = \langle X, z\rangle = X\langle
z\rangle$, and note that $P' \le T$ and so $T \trianglelefteq P$.
As $z$ does not centralize $P'= \langle x^2\rangle$,
$x^z \neq x^{1+2^{n-1}}$; hence (see for instance Corollary
 5.4.2 in \cite{go}) we have two possibilities:
 $x^z = x^{-1}$ and $T$ is dihedral, or
$x^z = x^{-1 + 2^{n-1}}$ and $T$ is semi-dihedral.

Let $u\in C$; then (since $u^2\in P'$),
$$\vert \langle u\rangle C'\vert \leq
\vert \langle u\rangle P'\vert = \frac {\vert u\vert\vert P'\vert}{\vert
\langle u\rangle \cap P'\vert}\leq 2\vert P'\vert = \vert X\vert.$$
This shows that $X/C'$ is a cyclic subgroup of maximal order in $C/C'$; so
there exists a complement $R/C'$ of $X/C'$ in $C/C'$. Then $XR = C$ and
$X\cap R = C'$.

Let $r\in R$, and suppose $r^2\neq 1$. Then $r^2\in R\cap P' = C' \leq Z =
\langle x^{2^{n-1}}\rangle$, and so $r^2 = x^{2^{n-1}}$. Recalling that
(as $n\geq 3$), $x^{2^{n-2}} \in P'$ is centralized by $r$, we have $\vert
rx^{2^{n-2}}\vert = 2$, and thus $r\in X\Omega_1(C)$.
Therefore, $C = X\Omega_1(C)$ and, consequently,
$ P = \langle C,z\rangle = \langle x,z\rangle\Omega_1(C) = T\Omega_1(C).$
Furthermore, $\vert \Omega_1(C)' \vert \le \vert C' \vert \le 2$,
and so $\Omega_1(C)$ has at most exponent $4$. In particular
$T \cap \Omega_1(C) = X \cap \Omega_1(C) \le \langle x^{2^{n-2}} \rangle \le Z(C)$.
Now let $g \in P$ and $y \in C$, with $\vert y \vert = 2$. Then $[g,y] \in T \cap \Omega_1(C)$,
so $[g,y,y]=1$ and hence $[g,y]^2 = [g,y^2] = 1$. But
$[P,y] \le P'$ is cyclic and so $[g,y] \in \Omega_1(P') = Z$.
Thus we have $[P,y] \le Z$ and so $[P, \Omega_1(C)] \le Z$; in particular every
subgroup of $\Omega_1(C)$ containing $Z$ is normal in $P$.

To finish, note that $\langle x^{2^{n-2}},
\Omega_1(C)\rangle/Z$ is a torus, so we may choose in it a
complement $D/Z$ of $\langle x^{2^{n-2}}\rangle/Z$, in such a way that
$\Omega_1(Z(P)) \leq D \leq \Omega_1(C)$. By what we observed above,
$D\trianglelefteq P$;
furthermore, $TD = T\Omega_1(C) = P$ and $T\cap D \leq T\cap \Omega_1(C) =
Z$. Since $Z \le T$ is central, has order 2 and 
$Z \le \Omega_1(Z(P)) \le D$, we have $Z \le T \cap D$ and hence $T \cap D = Z$.
Finally we observe that $D'\leq C' \leq Z$, has at most order $2$, whence $D = Z(D)E$
where $E=1$ or $E$ is extra-special. \ $\square$


\section{Upper intervals in $p$-groups}

An element $a$ of a poset $Q$
is named a \emph{conjunctive element} if, for all $x \in Q$,
there exists sup$\{a,x\}$ in $Q$. Quillen proved in
\cite{qu} that if $Q$ has a conjunctive element then
$\Delta(Q)$ is contractible.

Let $P$ be a $p$-group and $X\in \mathcal{A}_p(P)$. If
 $X \not \ge \Omega_1(Z(P))$, it is easily seen that
$X\Omega_1(Z(P))$ is a conjunctive element for
$\mathcal{A}_{p}(P)_{>X}$, and so we have that
$\Delta(\mathcal{A}_{p}(P)_{>X})$ is
contractible and therefore it is $(rk(P)-rk(X)-1)$-spherical.
Thus the only upper intervals $\mathcal{A}_p(P)_{>X}$ that we need to
study are those for $X\geq \Omega_1(Z(P))$.

Summarizing, for our purposes, we only have to consider
$\mathcal{A}_{p}(P)_{>X}$, for $P$ a $p$-group such that $\Omega_1(P)=P$,
$P'$ is cyclic and $X \ge \Omega_1(Z(P))$.

\n {\sc Upper intervals in extra-special groups.}
\ In this subsection, we
assume that $P$ is an extra-special $p$-group with center $Z= \langle
z\rangle\ ( = \Omega_1(Z))$. It is then well known that the
commutator $[\cdot,\cdot]$
induces a bilinear form $f:V \times V \to \mathbb{Z}_p$
on the $\mathbb{Z}_p$-vector space $V=P/Z$, defined by
$z^{f(u,v)}=[x,y]$, for $u=xZ$, $v=yZ$. With no loss of generality, we will also assume that
$P = \Omega_1(P)$ (thus, for odd $p$, $\exp (P) = p$).

\n We recall that if $P$ is such an extra-special group, and $Z< P_1 <
P$ is a subgroup of $P$ which is also extra-special, then, by letting
$P_2/Z$ be the orthogonal complement of $P_1/Z$ (with respect to the above
defined $f$), $P_2$ is also extra-special, and $P$ decomposes as the central
product $P= P_1\circ P_2$ (``$\circ$'' is our notation
for the central product).

We denote by
$D_8$ and $Q_8$, respectively, the dihedral and the quaternion group of
order $8$.


\begin{propos}\label{41}
Let $p$ be a prime and $P$ an extra-special $p$-group with center $Z=Z(P)$,
and let $X \in \mathcal{A}_p(P)_{\ge Z}$ such
that $X$ is not a maximal torus in $P$.
Then  $\mathcal{A}_{p}(P)_{>X}$ is $(rk(P)-rk(X)-1)$-spherical.
\end{propos}

\noindent \emph{Proof.} Clearly we have
$\mathcal{A}_p(P)_{>X} = \mathcal{A}_p(C_{P}(X))_{>X}$. If $P=C_P(X)$,
then $X=Z$ and  we can proceed as indicated in the points (1) and (2) below. 

Otherwise, we assume that $C_P(X)<P$ and use the bilinear form defined above
to reduce to the previous case. With respect to the bilinear form $f$ described above, $C_P(X)/Z$ is the
largest subspace of $P/Z$ orthogonal to $X/Z$, and since $X$ is not a maximal torus, we have $C_P(X)> X$.
Let $P_2/Z$ be a complement of $X/Z$ in $C_P(X)/Z$, hence we can rewrite $C_P(X)=XP_2$.
Clearly we have $P_2\neq P$, and, since $f$ is not
singular on $P/Z$, then $f$ is not singular on $P_2/Z$. Hence, $P_2$ is
extra-special, $Z(P_2)=Z$, and we may decompose $P$ as the central product $P_1\circ
P_2$, where $P_1/Z$ is the orthogonal complement of $P_2/Z$, and $P_1\cap P_2 =
Z$.

If $S \in \mathcal{A}_p(C_P(X))_{>X}$ is a maximal $p$-torus then, 
up to choosing another complement $\widetilde{P}_2/Z$ of $X/Z$ in $C_P(X)/Z$, we can assume that
$S=X K$, with $K \le P_2$ and $X \cap K=1$. Since $S$ is maximal in $C_P(X)$, the subgroup $KZ$ must be maximal in $P_2$,
which implies $rk(P_2)=rk(K)+1$. Thus,
$rk(P)=rk(S)=rk(X)+rk(K)=rk(X)+rk(P_2)-1$ and so $rk(P_2)=rk(P)-rk(X)+1$.

We now define the map:
$$
\begin{array}{crcl}
\varphi : & \mathcal{A}_p(P_2)_{>Z} & \rightarrow & \mathcal{A}_p(C_{P}(X))_{>X}  \\
	 & T & \mapsto & \varphi(T) = TX
\end{array}$$
By the Dedekind modular law it is easy to check that $\varphi$ is a poset
isomorphism and so
$\mathcal{A}_{p}(P)_{>X} \cong \mathcal{A}_{p}(P_{2})_{>Z}$ and we have reduced again
the study to an upper interval with respect to a center. We distinguish the cases $p > 2$ and  $p=2$.

\n \emph{(1) $p>2$.} We observe that 
$\mathcal{A}_{p}(P_{2})_{>Z}$ is the poset of inverse images of non-zero \emph{isotropic} subspaces $U$ of $V$
(they are the subspaces $U$ for which $f \vert_U = 0$).

\vspace{2mm}
\n \emph{(2) $p=2$.} \ In this case
we have
$P_2 \cong D_8 \circ \ldots \circ D_8 \circ D_8$ or
$P_2 \cong D_8 \circ \ldots \circ D_8 \circ Q_8$
(see e.g. Theorem 5.2 in \cite{go}). We observe that
$\mathcal{A}_2(P_2)_{>Z}$ is the set of the inverse images of the
isotropic and \emph{totally singular} subspaces of $V$ (they are the subspaces $U$ for which
the quadratic form $q(\alpha) = \alpha^{2}$ vanishes).

In both cases the complex $\mathcal{A}_{p}(P_{2})_{>Z}$ is known to be a building of rank $rk(P_2)-1$
(see example 10.4 in \cite{qu}). Thus, by the
Solomon-Tits Theorem (see \cite{brown2}, Theorem 5.2, p. 93), 
it is spherical of dimension $((rk(P_2)-1)-1)$. We observe that
$rk(P_2)-2=rk(P)-rk(X)+1-2$ and we are done. $\square$


Next, an easy extension of the previous result.
\begin{propos}\label{42}
Let $p$ be a prime, and $T$ a $p$-group, such that $\Omega_1(T)=T$,
$Z(T) < T$, $T/Z(T)$ is abelian and $T'$ is a cyclic group.
If $p=2$ we suppose further that $\vert T' \vert = 2$ and $Z(T)=\Omega_1(Z(T))$.
Then we have that
$\mathcal{A}_p(T)_{>\Omega_1(Z(T))}$ is
$(rk(T)-rk(\Omega_1(Z(T)))-1)$-spherical.
\end{propos}

\noindent \emph{Proof.} By assumption,
 $T' \le Z(T)$, and, since $\Omega_1(T)=T$,
$V=T/Z(T)$ is a $p$-torus.
If $p > 2$ then $T$ has exponent $p$, because it has class 2
(see e.g. Lemma 3.9 in \cite{go}); so $\vert T' \vert = p$ and it is generated by the elements of
order $p$. If $p=2$ we have by hypothesis that $\vert T' \vert = 2$.

We define a bilinear form $f:V \times V \to \mathbb{Z}_p$,
in the same way as for extra-special groups, and obtain again
that $T=Z(T)E$ with $E=1$ or $E$ extra-special. Thus
$T$ behaves as $C_P(X)$ in Proposition \ref{41}, and so
$\mathcal{A}_p(T)_{>\Omega_1(Z(T))} \cong \mathcal{A}_p(E)_{>T'}$.
By putting $X=Z(T)$ we can follow the same procedure of Proposition \ref{41},
where $E=P_2$, and we are done. $\square$

\n Let $C \circ D$ be a central product of a cyclic group $\langle a \rangle
= C \cong C_4$, and $D$ an  extra-special $2$-group of center
$\langle a^2 \rangle$.
We are going to study the interval
$\mathcal{A}_2(C\circ D)_{>Z}$, as we will later need
its homotopy type. Recall that $D \cong D_8 \circ \ldots \circ D_8 \circ D_8$
or $D \cong D_8 \circ \ldots \circ D_8 \circ Q_8$.

\n Let
$Ab(D) = \{ H \le D \mid H'=1 \}$ be the poset of all abelian subgroups
of $D$. We define
a map \ $
\pi_D: \mathcal{A}_{2}(CD)_{>Z} \to Ab(D)_{>Z}
$,
by letting $\pi_D(A)/Z$ be the projection of $\frac{A}{Z}$
on the second component of  $\frac{C}{Z} \times \frac{D}{Z}$.
This is an order preserving map and it has an inverse map $
\lambda: Ab(D)_{>Z} \to \mathcal{A}_{2}(CD)_{>Z}
$ which is defined by the rule
$$\lambda(X)=
\begin{cases}
X & \text{if $X$ is $2$-torus} \\
\langle ay, \Omega_1(X) \rangle & \text{if $X=\langle y, \Omega_1(X) \rangle$}
\end{cases}$$
and is also an order preserving map. Thus
$\mathcal{A}_{2}(CD)_{>Z} \cong Ab(D)_{>Z}$ and we are lead to
study the homotopy type of $Ab(D)_{>Z}$.

\begin{lemma}\label{43}
Let $D$ be a extra-special $2$-group of order $2^{2n + 1}$, and center $Z$.
Then $Ab(D)_{>Z}$ is $(n-1)$-spherical.
\end{lemma}

\noindent \emph{Proof.} Let $C= C_D(x)$ be the centralizer of
a non-central element $x$ of $D$. We define $Ab(C)_{>Z}$
to be an $Ab$-poset. Note that $Ab(C)_{>Z}$ is contractible because
it contains the conjunctive element
$\langle x \rangle Z$. We prove the following claim:

\emph{($\ast$) The union of any number
$m$ of $Ab$-posets in an extra-special
$2$-group of order $2^{2n + 1}$ is $(n-1)$-spherical.}

We will work on several induction arguments all of which are based on an
induction on the pair $k = (m,n)$ in lexicographic order.

For $j=1, \dots, m$, let $\mathcal{A}_j =
Ab(C_D(x_j))_{>Z}$ be $Ab$-posets of $D$. We write
$C_D(x_j) = C_j$. We observed above that each $\mathcal A_j$ is
contractible, so claim ($\ast$) is true for $m=1$.
We may then suppose  $m > 1$.

\emph{Claim 1.} $\mathcal{A}_i \cap \mathcal{A}_j$ is empty or 
$(n-2)$-spherical.

\emph{Proof of Claim 1.} Observe that $\mathcal{A}_i \cap \mathcal{A}_j = Ab(C_i \cap C_j)_{>Z}$.
If $[x_i,x_j]=1$, then $\langle x_i \rangle Z$ is a conjunctive element
in $\mathcal{A}_i \cap \mathcal{A}_j$, which is therefore contractible.
If $[x_i,x_j] \ne 1$, we let $D_\omega/Z$ be the orthogonal complement of
the subspace $\langle x_i,x_j \rangle/Z$ in $C_D(\langle x_1,x_2 \rangle)/Z$, 
with respect to the bilinear form defined at the beginning of the subsection. 
Thus $D_\omega = Z$ or $D_\omega$ is extra-special with
$rk(D_\omega)=rk(D)-1$, and $C_t=\langle x_t \rangle D_\omega$ for $t=i,j$.
It follows that $\mathcal{A}_i \cap \mathcal{A}_j = Ab(D_w)_{>Z}$
is empty or $(n-2)$-spherical by an inductive argument. $\square$

Let $m=2$. We have two possibilities:

(a) \ $[x_1,x_2]=1$;
then $\mathcal{A}_1$, $\mathcal{A}_2$ and $\mathcal{A}_1 \cap \mathcal{A}_2$
are all contractible, and so $\mathcal{A}_1 \cup \mathcal{A}_2$ is contractible.

(b) $[x_i,x_j] \ne 1$. If $D_\omega=Z$, then $rk(D)=2$ whence
$Ab(D)_{>Z}$ is a set of points, which is $0$-spherical.
If $D_\omega$ is extra-special then we conclude that
$\mathcal{A}_i \cup \mathcal{A}_j$
is $(n-1)$-spherical by the Gluing Lemma \ref{gl}.

Let now $m > 2$. If $[x_1,x_j]=1$, for all $j>1$, let $\mathcal U =
\mathcal{A}_2 \cup \dots \cup \mathcal A_m$; then
$\langle x_1 \rangle Z$ is a conjunctive element in
$ \mathcal{A}_1 \cap \mathcal U$, which is therefore contractible,
and $\mathcal U$ is $(n-1)$-spherical by
induction. By Lemma \ref{gl}, $\mathcal{A}_1 \cup \mathcal U$
is $(n-1)$-spherical. 

Otherwise, we can suppose $[x_1,x_2] \ne 1$.
For each $j=2, \ldots, m$ we put $\mathcal{B}_j:=\mathcal{A}_1 \cap \mathcal{A}_j$,
and set $\mathcal{U} := \mathcal{A}_3 \cup \dots \cup \mathcal A_m$.
By the initial remark and an inductive argument we have that both
$\mathcal{A}_1$ and $\mathcal{A}_2 \cup \mathcal{U}$ are
$(n-1)$-spherical. So we have to consider
$\mathcal{A}_1 \cap (\mathcal{A}_2 \cup \mathcal{U})=
 \mathcal{B}_2 \cup \dots \cup \mathcal{B}_m$.

\emph{Claim 2:} $\mathcal{A}_1 \cap (\mathcal{A}_2 \cup \mathcal{U})$ is empty or $(n-2)$-spherical.

\emph{Proof of Claim 2.} We want to prove that any union  $\mathcal{B}_{r_1} \cup \dots \cup \mathcal{B}_{r_v}$
with $r_j \in \{2,\ldots,m \}$ is either empty or $(n-2)$-spherical and we use induction on the number of $\mathcal{B}_{r_j}$'s.
We consider a union $\underset{j=1}{\overset{t} {\bigcup}} \mathcal{B}_{r_j}$, for some $1 \le t \le m-1$.
If $[x_1,x_{r_j}]=1$ for all $j=1, \ldots, t$, then $\langle x_1 \rangle Z$ is a conjunctive element for 
$\underset{j=1}{\overset{t} {\bigcup}} \mathcal{B}_{r_j}$, which is thus contractible. Otherwise,
there is a $j_0$ such that $[x_1,x_{r_{j_0}}] \ne 1$.

Without any loss of generality and up to relabeling, we can assume that $r_{j_0}=2$
and the union is $\underset{j=2}{\overset{t} {\bigcup}} \mathcal{B}_j$, for some $2 \le t \le m$.
By Claim 1, each $\mathcal{B}_j$ is either empty or $(n-2)$-spherical and by induction hypothesis on the number of
the $\mathcal{B}_j$'s, the same holds for 
$\underset{j=3}{\overset{t} {\bigcup}} \mathcal{B}_j$. 
We set $\mathcal{U}_t := \mathcal{A}_3 \cup \dots \cup \mathcal A_t$ and we study the homotopy type of
$\mathcal{B}_2 \cap \Big(\underset{j=3}{\overset{t} {\bigcup}}
\mathcal{B}_j \Big) =
\mathcal{A}_1 \cap \mathcal{A}_2 \cap \mathcal{U}_t$.

\n We want to classify the structure of
$\mathcal{A}_1 \cap \mathcal{A}_2 \cap \mathcal{A}_j$, for
$3 \le j \le t$. If there is a $3 \le j \le t$ such that
$x_j \in \langle x_1,x_2 \rangle$, then
$C_1 \cap C_2 \cap C_j = D_\omega$ and so
$\mathcal{A}_1 \cap \mathcal{A}_2 \cap \mathcal{A}_j =
\mathcal{A}_1 \cap \mathcal{A}_2$ and
$\mathcal{A}_1 \cap \mathcal{A}_2 \cap \mathcal{U}_t =
\mathcal{A}_1 \cap \mathcal{A}_2$, this implies that
$(\mathcal{A}_1 \cap \mathcal{A}_2) \cup (\mathcal{A}_1 \cap \mathcal{U}_t)
= \underset{j=2}{\overset{t} {\bigcup}} \mathcal{B}_j$
is either empty or $(n-2)$-spherical and we are done.
Otherwise, for all $3 \le j \le t$, we have that
$x_j \notin \langle x_1,x_2 \rangle$, hence
$x_j = r_j s_j$, with $r_j \in \langle x_1,x_2 \rangle$
and $s_j \in D_\omega \backslash Z$.
Then $C_1 \cap C_2 \cap C_j =
C_{D_\omega}(s_j)$, and so $\mathcal{A}_1 \cap \mathcal{A}_2 \cap \mathcal{U}_t$
can be written as a union that has the same structure of
$\underset{j=1}{\overset{m} {\bigcup}} \mathcal{A}_j$. Hence we have recovered
the same type of union of the thesis of the Lemma on a smaller scale. Since
$rk(D_\omega) < rk(D)$ we conclude, by induction on $k$,
that $\mathcal{A}_1 \cap \mathcal{A}_2 \cap \mathcal{U}_t$ is either
empty or $(n-3)$-spherical.

Now, using backward the Gluing Lemma, we have that
$\mathcal{A}_1 \cap (\mathcal{A}_2 \cup \mathcal{U}_t)$ is either
empty or $(n-2)$-spherical. Since this is true for any union, Claim 2 is now proved. 
$\square$

To finish the proof of Claim ($\ast$) we observe that Claim 2 and the Gluing Lemma imply that
$\underset{j=1}{\overset{m} {\bigcup}} \mathcal{A}_j$ is
$(n-1)$-spherical. Finally we observe that
$Ab(D)_{>Z}$ is a union of $Ab$-posets since, if $H \le D$,
$H'=1$ and $x \in H$, then $H \le C_D(x)$. By  ($\ast$)
the proof is now complete.
$\square$

\begin{corol}\label{44}
Let $P = C \circ D$ be a central product of a cyclic group $C$ of order 4
and an  extra-special $2$-group $D$ with center $Z = Z(D)$. 
Let $X \in \mathcal{A}_p(P)_{\ge Z}$ such
that $X$ is not a maximal torus in $P$.
Then  $\mathcal{A}_{p}(P)_{>X}$ is $(rk(P)-rk(X)-1)$-spherical.
\end{corol}

\begin{osserv}\label{45}
\emph{By following steps similar to those of Claim ($\ast$)
of the previous Lemma, it is possible to directly prove that:
\emph{if $P=\Omega_1(P)$
is an extra-special $p$-group, with $p$ any prime number, then
$\mathcal{A}_p(P)_{>Z(P)}$ is $(rk(P)-2)$-spherical}. The argument makes no use of the Solomon-Tits
Theorem on buildings. We omit the details of this proof.}
\end{osserv}

\n {\sc The general case.} \ Now we describe upper intervals $\mathcal
A_p(P)_{>X}$ for those $p$-groups $P$ we are interested in. The results in
section $3$, at least for $p$ odd, immediately reduce to the extra-special
case.
\begin{propos}\label{51}
Let $p$ be an odd prime number and $P$ a $p$-group such that
$\Omega_1(P)=P$, $P' \ne 1$ is cyclic and set $X = \Omega_1(Z(P))$.
Then $\mathcal{A}_p(P)_{>X}$ is $(rk(P)-rk(X)-1)$-spherical.
\end{propos}

\noindent \emph{Proof.} By Theorem \ref{31}, we have that $\vert P' \vert = p$,
so $P'$ is cyclic and $P' \le Z(P)$. Thus $P/Z(P)$ is abelian. Now
Proposition \ref{42} completes the proof. $\square$

\n In the case $p=2$ the situation is more complicated and we restrict ourselves to
the case of Theorem \ref{32} with a central product. We recall from
Theorem \ref{32} that if $p=2$ then, under suitable hypotheses, 
$P=TD$ where $T$ is dihedral or semi-dihedral
and $D$ is a product of a $2$-torus and an extra-special group.

\begin{propos}\label{52}
Let $P$ be a $2$-group such that
$\Omega_1(P)=P$ and $P' \ne 1$ is cyclic. Set $X = \Omega_1(Z(P))$,
and if $|P'|>1$ let $P=TD$ be the decomposition of Theorem \ref{32}. Then
\begin{enumerate}
\item If one of the following holds:
\vspace{1mm}

$(i)$ $\vert P' \vert =2$ and $X=Z(P)$, or

\vspace{1mm}
$(ii)$ $\vert P' \vert >2$ and $T$ is dihedral,

\vspace{1mm}
\n then $\mathcal{A}_2(P)_{>X}$ is $(rk(P)-rk(X)-1)$-spherical.

\item if $\vert P' \vert >2$, and $T$ is semidihedral,
$\mathcal{A}_2(P)_{>X}$ is not weakly spherical.
\end{enumerate}
\end{propos}

\noindent \emph{Proof.} If $\vert P'\vert = 2$, then the conclusion follows via
Proposition \ref{42} as in the odd case.

Thus, let $\vert P' \vert > 2$ and define
$C=C_P(P')$.
Let $\langle x \rangle = F$ be a maximal cyclic subgroup of $T$, with
$F\trianglelefteq P$; if $\vert x \vert = 2^n$, let $Z = Z(T) = \langle
x^{2^{n-1}}\rangle$.
In both cases (see e.g. Theorem 4.3 in \cite{go})
$T/Z$ is a dihedral group.

Let $A$ be a
torus of $P$, then $(AD \cap T)' \le (AD)' \cap T \le Z(D) \cap T \le Z(T)=Z$
so $(AD \cap T)/Z$ is abelian and therefore $(AD \cap T)/Z$ is a torus
of $T/Z$. Thus $AD \cap T$ is contained in a suitable
$R= \langle x^{2^{n-2}}, t \rangle$ of order $8$,
for some $t \in T \backslash F$. 
Therefore we can say that
\emph{if $A\in \mathcal{A}_2(P)$ then
$A \in \mathcal{A}_2(RD)$ for a suitable $R \le T$ with
$\langle x^{2^{n-2}} \rangle \le R$, $R \cong D_8$ or $R \cong Q_8$}.


Suppose now that $T$ is dihedral. Then (see \cite{as}) no
subgroup of $T$ is isomorphic to $Q_8$. We take
$\mathcal{D} = \{ R \le T \mid \langle x^{2^{n-2}} \rangle \le R,
\; R \cong D_8 \}$. Since $\mathcal D$ is finite, elements
$R_1, \ldots, R_k \in \mathcal{D}$ exist such that
$$\Delta(\mathcal{A}_2(P)_{>X})=
\Delta(\mathcal{A}_2(R_{1}D)_{>X}) \cup \ldots \cup
\Delta(\mathcal{A}_2(R_{k}D)_{>X})$$
and
$R_{i}D \ne R_{j}D$ for every $1 \le i,j \le k$, $i \ne j$. We see that
$(R_{1}D) \cap (R_{2}D) = (R_1 \cap R_2)D= \langle x^{2^{n-2}} \rangle D$.
Define $\Sigma_i := \Delta(\mathcal{A}_2(R_{i}D)_{>X})$
and $\Sigma:= \underset{j=1} {\overset{k}{\bigcup}} \Sigma_j$.
By induction on $k$, we prove that $\Sigma$ is $(rk(R_{1}D)-rk(X)-1)$-spherical.
If $k=1$, we have already proved it, since $T=R_1$ and so $|P'| \le 2$ and Proposition 4.2
completes the proof. If $k>1$ we put $\Sigma_0=
\underset{j=1} {\overset{k-1}{\bigcup}} \Sigma_j$ and,
by inductive hypothesis, we know that $\Sigma_0$ is
$(rk(R_{1}D)-rk(X)-1)$-spherical.
Thus we observe that $\Sigma_0 \cap \Sigma_k =
\mathcal{A}_2(\langle x^{2^{n-2}} \rangle D)_{>X}$. If $Z(D)\setminus
\langle x^{2^{n-2}} \rangle$ contains
an element $y$ of order $4$,
we have that $\langle yx^{2^{n-2}}, x^{2^{n-1}} \rangle$ is a conjunctive
element for $\Sigma_0 \cap \Sigma_k$, which is therefore contractible. Otherwise,
we may have the following two cases: $Z(D) >Z$, which implies that $Z(D)$
is a conjunctive element, or $Z(D)=Z$ and so $D$ is extra-special and we can apply
Corollary \ref{44} to obtain that $\Sigma_0 \cap \Sigma_k$ is $(rk(D)-rk(X)-1)$-spherical.
We observe that $rk(D)-rk(X)-1=rk(R_1 D)-rk(X)-2$ so, by the Gluing Lemma, we have that
$\Sigma = \Delta(\mathcal{A}_2(P)_{>X})$ is $(rk(R_{1}D)-rk(X)-1)$-spherical.
Since $rk(R_1D)=rk(P)$, we can say that
\emph{if $T$ is dihedral,
then $\Delta(\mathcal{A}_2(P)_{>X}))$
is $(rk(P)-rk(X)-1)$-spherical.}

Suppose now that $T$ is semidihedral, we then know (see \cite{as})
that $T_1,T_2 \le T$ exist such that
$T_1 \cong D_8$ and $T_2 \cong Q_8$.
It can happen
that $rk(T_{1}D)=rk(P)$ or $rk(T_{2}D)=rk(P)$. We will work
out the case $rk(T_{1}D)=rk(P):=r+1$ and $rk(T_{2}D)=r$, the other being
similar.  Let $A_1$ be a maximal torus of $T_{1}D$,
and $A_2$ a maximal torus of $T_{2}D$ (so that
$rk(A_2)=rk(P)-1$).
Then, by Theorem 12.1 in
\cite{qu}, applied on $A_1$ and $A_2$,
we have that $\widetilde{H}_r(\mathcal{A}_2(S)) \ne 0 \ne
\widetilde{H}_{r-1}(\mathcal{A}_2(S))$, and so $\mathcal{A}_2(S)$
is not even weakly spherical.
Hence we have that \emph{if $T$ is semidihedral, then $\Delta(\mathcal{A}_2(P))$ is not
weakly spherical.} $\square$


\section{Homotopy type of the Quillen complex}

In this section, we discuss the sphericity of the whole Quillen complex
$\mathcal A_p(G)$, for solvable groups $G$ whose Sylow $p$-subgroups have a
cyclic derived group. But first, we need to get some information about the way
the Sylow $p$-subgroups are located in such groups.

Let $G$ be a group and $p$ a prime. We denote by $O_p(G)$ the largest normal $p$-subgroup of $G$ (it may be the trivial
group). Similarly, we denote by $O_{p'}(G)$ the largest normal subgroup of $G$ whose order is not divisible by $p$.
We denote by $O_{p',p}(G)$ the inverse image in $G$ of the group $O_p(G/O_{p'}(G))$ and, in a
similar fashion, we define the subgroups $O_{p',p,p'}(G),O_{p',p,p',p}(G)$ and so on, building an ascending normal series
(called the $p$-\emph{series} of $G$). If this ascending series reaches $G$, we say that $G$ is $p$-\emph{solvable} and 
we define $\ell_p(G)$ (the $p$-\emph{length} of $G$) to be the number factors of the $p$-series that are $p$-groups. 
We now need to find an upper bound for $\ell_p(G)$ when the Sylow
$p$-subgroups of $G$ have cyclic derived groups
(if the Sylow $p$-subgroups $P$ of $G$
are abelian, then $\ell_p(G)$ is easily seen to be $1$). For a group $G$ we define $Fit(G)$
(the \emph{Fitting subgroup} of $G$) to be the largest normal nilpotent subgroup of $G$. 
For $p$-solvable groups it might happen that $Fit(G)=1$, so it is possible to define the so-called
\emph{generalized Fitting subgroup} for a $p$-solvable group $G$ which is always non-trivial. 
However, we will work with a $p$-solvable group $G$ such that $O_{p'}(G)=1$ and in this case 
the generalized Fitting subgroup
coincides with $Fit(G)$. The propositions we state
are well known (for instance as applications of the Hall-Higman Theorem).
Anyway, for the sake of completeness, we will directly prove them.

\begin{lemma}\label{LP1}
Let $G$ be a finite $p$-solvable group, for a prime number $p \ge 5$.
Let $P \in \text{\emph{Syl}}_p(G)$
and suppose that $P'$ is cyclic. Then $\ell_p(G) = 1$.
\end{lemma}

\emph{Proof.} With no loss of generality, we can assume $P' \ne 1$ and
$O_{p'}(G)=1$, so that $N=O_p(G) = Fit(G)$. Since $\Phi(G) \le Fit(G)$
we may also assume that $\Phi(G)=1$. Furthermore, we notice that $Fit \left( G/\Phi(G) \right)
=Fit(G)/\Phi(G)$ (see 5.2.15 of \cite{rob}), so we may suppose that $N$ is a $p$-torus.

We argue by contradiction that $N < P$ (i.e. $P$ is not normal in $G$), and
choose an element $g \notin N_G(P)$. Let $M= \langle P' \cap N,
(P')^g \cap N \rangle$. If $h \in G$, then $1 \ne [P^h,N] \le (P')^h$,
so that $[P^h,N]=(P')^h \cap N \cong C_p$ and $M$
is abelian of rank at most $2$. We observe that
$P' \cap N \le M \le N \unlhd P$ and $MP' \unlhd P$, thus
by Dedekind modular law we have $M=M(P' \cap N)=MP' \cap N \unlhd P$. Similarly,
we can prove that $M \unlhd P^g$, whence $M$ is
normalized by $S= \langle P, P^g \rangle$.

Let $C=C_S(M)$, and $H \le C$ a Hall $p'$-subgroup of $C$. We observe
$[N,H] = [N,H,H] \le [N,C,C] \le [[N,S],C] \le [M,C] = 1$, so
$H \le C_G(N) = N$ which is a $p$-group. Thus, $H=1$ and $C$ is a
$p$-group. Since $M$ is normal in $S$, $C\unlhd S$ and so $C\leq P\cap P^g$.

If $\vert M \vert = p$, the quotient group 
$S/C$ embeds in $Aut_{\mathbb{F}_p}(M)$ which is cyclic of
order $(p-1)$. Since $S$ is generated by Sylow $p$-subgroups, this forces
$S=C$, and so we get the contradiction $P = C = P^g$. We point out that
the proof holds for any prime $p$ up to here.

Hence we must have $\vert M \vert = p^2$, and so $S/C$ embeds in $Aut_{\mathbb{F}_p}(M) = GL(2,p)$.
Arguing as in the previous case, we have that $p$ divides $\vert S/C
\vert$. So $\vert P/C \vert = p$ and
$O_p(S/C) = 1$ (otherwise $P=P^g$).
Since it is generated by its Sylow $p$-subgroups, the group $S/C$ embeds in $SL(2,p)$:
checking the list of $p$-solvable subgroups of $SL(2,p)$ (see e.g.
6.25 and 6.26 in \cite{su1}) we find that $O_p(S/C) \ne 1$ and this is
a contradiction.

This shows that $P \trianglelefteq G$ and so $\ell_p(G)=1$. $\square$

The previous Lemma does not hold for $p < 5$. For instance, the group $S_4$ has
$2$-length 2 and Sylow 2-subgroups isomorphic to $D_8$; another example is the group
$(C_3 \times C_3) \rtimes SL(2,3)$ that has $3$-length 2 and
Sylow 3-subgroups extra-special of order 27 and exponent 3.
In a sense, these are the only cases that can occur.

\begin{lemma}\label{LP2}
Let $p$ be a prime, $p < 5$. Let $G$ be a finite solvable group such
that $p$ divides its order. Let $P \in \text{\emph{Syl}}_p(G)$ and
suppose that $P'$ is cyclic. Then $\ell_p(G) \le 2$.
If $\ell_p(G) = 2$, we assume $O_{p'}(G)=1$.
Let $N = O_{p}(G)$ and $g \in G$ be a fixed element
such that $g \notin N_G(P)$. Then:

$(i)$ $\frac{\langle P, P^g \rangle}{N} \cong SL(2,2)$, if $p=2$.

$(ii)$ $\frac{\langle P, P^g \rangle}{N} \cong SL(2,3)$, if $p = 3$.
\end{lemma}

\emph{Proof.} We argue as in the proof of the previous Lemma and
assume that $O_{p'}(G)=1$ and let $N= O_p(G) = Fit(G)$. Again, as in the previous
proof, we can assume that $\Phi(G)=1$, so that $N$ is a $p$-torus.

If $\ell_p(G)=1$ there is nothing to prove. Thus, let $\ell_p(G) > 1$
(and therefore $P' \ne 1$) and set $S= \langle P,P^g \rangle$. Let
$M = \langle P' \cap N, (P^g)' \cap N \rangle$, and $C=C_s(M)$. 

Arguing as in the proof of the previous Lemma, we have
$\vert M \vert = p^2$ and $C$ is a normal $p$-subgroup of $S$.

If $[M,P]=1$, then $P \le C$ and so $P = C = S$ and we get a contradiction
because $P=P^g$. Hence $[M,P] \ne 1 \ne [M,P^g]$ and, as in Lemma \ref{LP1}, $[P^h,N]= (P^h)' \cap N \cong C_p$,
for all $h \in G$. It follows
that $M \le [S,M] \le [S,N] \le M$ and so $M=[N,S]$.

Now $S/C \le Aut(M)=GL(2,p)$ and since it is generated by its
Sylow $p$-subgroups, we must have $S/C \le SL(2,p)$. Since $S/C$
admits two distinct Sylow $p$-subgroups, we have
\medskip \\
(i) $S/C \cong SL(2,2)=S_3$, if $p=2$
\medskip \\
(ii) $S/C \cong SL(2,3)$, if $p=3$.
\medskip \\
Consider $T/C = O_{p'}(S/C)$, which is a non-trivial group, and let $T_0$ be a 
$p'$-subgroup of $T$ such that $T=C \cdot T_0$. We see that $[N,P] \le N \cap P' =
\Omega_1(P)$, the last subgroup being cyclic of order $p$ and contained in $Z(P)$. Hence we have
$[N,P,P]=1$. By the Three Subgroups Lemma we have $[P',N]=1$. Since $N=Fit(G)=C_G(N)$,
it follows that $P' \le N \cap P' \le M$. Thus $[C,P] \le P' \le M$,
$[C,P^g] \le M$ and $[C,S] \le M$. This implies
$[C,T_0] \le M \le T_0  M$ and so 
\[
\frac{T}{M}=\frac{C}{M} \times \frac{T_0 M}{M}.
\]
We observe that $T_0 M / M$ char $T/M \trianglelefteq S/M$ and so
$T_0 M \trianglelefteq S$. Let $Y=C_N(T_0)$, then
$Y=C_N(T_0 M) \trianglelefteq S$. By Theorem 2.3 on page 177
in \cite{go}, we have that
\[
N= Y \times [N,T_0]
\]
with $[N,T_0] \le M$
and $1 \ne [N,T_0]=[N,T_0 M] \trianglelefteq S$. Thus $M= [N,T_0]$ and 
so $N= Y \times M$, with $Y,M \trianglelefteq S$. Now we observe that 
\[
[N,C] = [YM,C] = [Y,C] \le [N,S] \cap Y = M \cap Y = 1
\]
and so $C \le C_G(N) = N$. By definition of $N$ it is clear that $N \le C$
and so it follows that $N=C$.

If $P$ is a Sylow $p$-subgroup of $G$ then, by  statements (i) and (ii), we have $NP/N \cong C_p$ and so
$O_{p',p}\left(\frac{G}{N}\right) \le Z \left(\frac{PN}{N}\right) = 
\frac{PN}{N} \cong C_p$.
Thus $\ell_p(G)=2$. $\square$

Some easier comments are needed before returning to the Quillen complex.

\begin{lemma}\label{LP3}
Let $p < 5$ be a prime number and $G$ a finite solvable group such that
$\ell_p(G) = 2$. Let $P \in \text{\emph{Syl}}_p(G)$ and suppose that $P'$
is cyclic. If $H= O_{p'}(G)$, $N=O_{p',p}(G)$ and $g \notin N_G(PH)$,
then $N = (PH \cap P^{g}H)$. Moreover, $N=HQ$ with $Q \le P$ and $[P:Q]=p$.
\end{lemma}
\noindent \emph{Proof.} Let $g$ be as in the statement; then
$gH \notin N_{G/H}(PH/H)$. Thus the hypotheses of Lemma \ref{LP2}
hold for $\langle P,P^g \rangle /N$.
This means that $[PH:N]=[P^{g}H:N]=p$ and, recalling that
$N \le (PH \cap P^{g}H) < PH$, the Lemma follows. $\square$

\begin{lemma}\label{LP4}
Let $p$ be a prime number, $P$ a $p$-group and $Q \le P$ such that
$[P:Q] = p$. Then $rk(P)-1 \le rk(Q) \le rk(P)$.
\end{lemma}
\noindent \emph{Proof.} Obvious.  $\square$


\begin{teor}\label{55}
Let $p$ be a prime number and let $G$ be a finite solvable group such that
$p$ divides its order. Let $P \in \text{\emph{Syl}}_p(G)$ and suppose that
$P'$ is cyclic. Suppose that $p>2$ or that $p=2, \Omega_1(P)=T \circ D$ as in Theorem \ref{32} with 
$T=1$ or $T$ is dihedral. Then $\mathcal{A}_p(G)$ is $(rk(P)-1)$-spherical.
\end{teor}

\noindent \emph{Proof.} We define $S=O^{p'}(G)$ and observe that $\mathcal{A}_p(G)=\mathcal{A}_p(S)$.
We need only to study $\mathcal{A}_p(S)$.

Suppose first that $\ell_p(S)=1$. Then, since $O^{p'}(S)=S$, we have
$S =O_{p',p}(S)$ and so $S$ can be written as $S= N \rtimes P$,
where $N= O_{p'}(S)$. We apply Propositions \ref{51} and \ref{52} and Lemma \ref{23} to conclude that
$\mathcal{A}_p(S)$ is $(rk(P)-1)$-spherical.

Suppose now that $\ell_p(S)>1$. Then, by Lemmas \ref{LP1} and \ref{LP2},
we have $p<5$ and $\ell_p(S)=2$. Let $P_1, \ldots,
P_k$ be all the distinct Sylow $p$-subgroups of $G$, and
set $T_i=O_{p'}(S)P_i$.
It may happen that $T_i = T_j$, for some $1 \le i,j \le k$, $i \ne j$;
thus, let us select a set  $T_{j_1}, \ldots, T_{j_r}$ of distinct
representatives of the $T_i$'s.
In order to keep the argument simple we rename $T_{j_v}=T_v$.

Then
$\mathcal{A}_p(S)= \mathcal{A}_p(T_1) \cup \ldots \cup \mathcal{A}_p(T_r)$.
Moreover, since each $T_i$ has $p$-length equal to $1$,
$\mathcal{A}_p(T_i)$ is $(rk(P)-1)$-spherical. By Lemma \ref{LP3} we have that
$T_i \cap T_j = O_{p',p}(S)$ and so 
\[
\mathcal{A}_p(T_i) \cap \mathcal{A}_p(T_j) =
\mathcal{A}_p(T_i \cap T_j) =
\mathcal{A}_p(O_{p',p}(S)) = \mathcal{A}_p(O_{p'}(S)Q)
\] 
for some suitable $p$-subgroup $Q \le P_i$ such that $[P_i:Q]=p$ .
By Lemmas \ref{23} and \ref{LP4}, we have that 
$\mathcal{A}_p(O_{p'}(S)Q)$ is $(rk(P)-1)$-spherical or $(rk(P)-2)$-spherical.
Applying the Gluing Lemma \ref{gl} to this
covering of $\mathcal{A}_p(S)$ completes the proof of the Theorem.
$\square$


\section{Homotopy type of the links of $\mathcal{A}_p(G)$}

In this last section we deal with the topological structure of
the links of the Quillen complex $\mathcal{A}_p(G)$, and prove that, in the cases stated in
Theorems \ref{main} and \ref{example} the Quillen complex is Cohen--Macaulay.
For this section we will assume that $p$ is a prime number and $G$ is a finite solvable group such that
$p$ divides its order. We let $P \in \text{Syl}_p(G)$ and suppose that
$P'$ is cyclic. Finally, we assume 
that $p>2$ or that $p=2, \Omega_1(P)=T \circ D$ as in Theorem \ref{32} with 
$T=1$ or $T$ is dihedral.

Recall that if
$\Gamma$ is a complex and $\sigma \in \Gamma$ then the link of $\sigma$ is
defined by
$L_\Gamma(\sigma) = \{ \tau \in \Gamma \mid \tau \cap
\sigma = \varnothing, \tau \cup \sigma \in \Gamma \}$.

If $\Gamma = \mathcal{A}_{p}(G)$, $\sigma=\{H_1 < H_2 < \ldots < H_k \}
\in \Gamma$, and we set $H_0=1$, then $L_\Gamma(\sigma)=
(H_0,H_1)\ast (H_1,H_2) \ast \ldots \ast \mathcal{A}_{p}(G)_{> H_k}$.
It is obvious that $(H_0,H_1) = \mathcal{A}_{p}(G)_{<H_1}$ and
that every $(H_i,H_{i+1})$ is $(rk(H_{i+1})-rk(H_i)-2)$-spherical,
since $(H_i,H_{i+1}) = \mathcal{A}_{p}(H_{i+1})_{>H_i} \setminus
\{ H_{i+1} \} \cong
\mathcal{A}_p(H_{i+1}/H_i) \setminus
\{H_{i+1}/H_i\}$ which is spherical (see the discussion at
the end of page 118 in \cite{qu} and Theorem 10.6  of
\cite{hc}).
Thus, in order to prove that
$L_\Gamma(\sigma)$ is spherical, we need to prove that
every $\mathcal{A}_p(G)_{>X}$ is $(rk(G)-rk(X)-1)$-spherical, for any $X \in \mathcal{A}_p(G)$.

The technical tool that we will need is the following 
criterion of homotopy equivalence, due to Pulkus and Welker (Corollary 2.4 in \cite{pw}).
\begin{teor}\label{pulwe}
Let $f : \mathcal P \rightarrow \mathcal Q$ be an order--preserving map of finite posets. Assume that

\vspace{1mm}
$(i)$ $\mathcal Q$ is a meet--semilattice, with $\hat 0 = \min \mathcal Q$;

\vspace{1mm}
$(ii)$ all elements $q\in \mathcal Q$, with the possible exception of $q= \hat 0$, belong to the image of $f$;

\vspace{1mm}
$(iii)$ for every $\hat 0 \ne q\in \mathcal Q$ there exists $c_q\in f^{-1}(\mathcal Q_{\le q})$ such that the inclusion map $\Delta (f^{-1}(\mathcal Q_{< q})) \rightarrow \Delta (f^{-1}(\mathcal Q_{\le q}))$ is homotopy equivalent to the constant map $c_q$.

\vspace{1mm}
Then 
$$
\Delta (\mathcal P) \simeq
\bigvee_{q\in \mathcal Q}\Delta (f^{-1}(\mathcal Q_{\le q})) \ast \Delta(\mathcal Q_{>q}).$$
\end{teor}

\noindent Let now $X \in \mathcal{A}_p(G)$.
We observe that $\mathcal{A}_p(G)_{>X} =
\mathcal{A}_p(C_G(X))_{>X}$, and so we may assume that
$X \le Z(G)$. We can suppose $X = \Omega_1(Z(G))$, otherwise $X$ would be a conjunctive element.
Moreover since $\mathcal{A}_p(G)_{>X} = \mathcal{A}_p(O^{p'}(G))_{>X}$,
we can also assume that $G=O^{p'}(G)$.

\begin{propos}
Suppose $X= \Omega_1(Z(G))$. Under the assumptions made on $G$, we have that
$\mathcal{A}_p(G)_{>X}$ is $(rk(G)-rk(X)-1)$-spherical.
\end{propos}

\noindent \emph{Proof.} The proof is  by induction on $\vert G \vert$. 
The idea is to first establish a sort of 
``local'' version of Pulkus and Welker's wedge decomposition (Theorem \ref{PW}), 
to which the inductive assumption is  then applied.
As some of the arguments are essentially the same as those used by 
Pulkus and Welker in \cite{pw}, we will be rather sketchy at some points.

Thus, let $N=O_{p'}(G) \ne 1$, and write
$\overline{G} = G/N$, and $\overline{A}=AN/N$ for each
$A \in \mathcal{A}_p(G)$. Let $\mathcal{P}=\mathcal{A}_p(G)_{>X}$, and $\mathcal{Q}=\mathcal{A}_p(\overline{G})_{>\overline{X}}
\cup \{ \overline{X} \}$. Then $\mathcal Q$ is a meet--semilattice with least element $\hat{0}=\overline{X}$.

Then, setting $f(A) = \overline{A}$ for all $A\in \mathcal{P}$, defines an order preserving map $f: \mathcal{P} \to \mathcal{Q}$.
Now, since $X=\Omega_1(Z(G))$ is a normal elementary abelian $p$-subgroup of $G$, and $N$ is a normal $p'$-group, we have that, for every $p$-subgroup $K$ of $G$, $KN \ge XN$ if and only if $K\ge X$; thus,
$$f^{-1}(\mathcal{Q}_{\le \overline{A}})= \mathcal{A}_p(AN)_{>X}
\cong \mathcal{A}_p(AN/X).$$
for all $A\in \mathcal P$. It then follows from Quillen's Theorem \ref{tq} that  $f^{-1}(\mathcal{Q}_{\le \overline{A}})$ is $(rk(A)-rk(X)-1)$-spherical, for all $A\in \mathcal P$.

In particular, $f^{-1}(\mathcal{Q}_{\le \overline{A}})$ is $(rk(A)-rk(X)-2)$-connected. Since, moreover, $f^{-1}(\mathcal{Q}_{< \overline{A}})$ is $(rk(A)-rk(X)-2)$-dimensional, standard topological arguments (see Lemma 3.2 in \cite{pw}, which is essentially what we need) entail that, for every $\overline A\in \mathcal Q \setminus \{\hat 0\}$, the inclusion map 
$$f^{-1}(\mathcal{Q}_{< \overline{A}})
\to f^{-1}(\mathcal{Q}_{\le \overline{A}})$$
is homotopic to a constant. We are then in a position to apply Proposition \ref{pulwe}, to obtain the wedge decomposition formula
$$
\mathcal{A}_p(G)_{>X} \simeq
\mathcal{A}_{p}(\overline{G})_{> \overline{X}} \vee
\underset{\overline{A} \in \mathcal{A}_{p}(\overline{G})_{> \overline{X}}}
\bigvee
\mathcal{A}_{p}\left(\frac{AN}{X}\right) \ast
\mathcal{A}_{p}(\overline{G})_{>\overline{A}}.$$
(where, clearly, we have extracted from the big wedge the contribution of $\hat 0 = \overline X$, i.e. 
$\mathcal{A}_{p}(\overline{G})_{> \overline{X}}$).

Now, by the inductive hypothesis on $\vert G \vert$,  we have that
$\mathcal{A}_{p}(\overline{G})_{> \overline{A}}$ is
$(rk(\overline{G})-rk(\overline{A})-1)$-spherical for each $\overline A \in \mathcal Q$.
If $\overline X < \overline A \in \mathcal Q$, then $\mathcal{A}_{p}\left(\frac{AN}{X}\right)$ is $(rk(A)-rk(X)-1)$--spherical by Quillen's Theorem, and so (since, clearly, $rk(\overline{G})-rk(\overline{A}) = rk (G) - rk (A)$),
$$\mathcal{A}_{p} (AN/X) \ast
\mathcal{A}_{p}(\overline{G})_{>\overline{A}}$$
 is $(rk (G)-rk(A)-1)$--spherical, for every $\overline X < \overline A \in \mathcal Q$. 
Hence, by the above formula, we conclude that $\mathcal{A}_{p}(G)_{> X}$ is
$(rk(\overline{G})-rk(\overline{X})-1)$-spherical, and case $N \ne 1$ is done.

Thus, assume now $N=1$, and recall that $\ell_p(G) \le 2$. 

If $\ell_p(G)=1$, this means that
$G$ itself is a $p$-group and so we are done by the results in section $4$
(namely, Propositions \ref{51} and \ref{52}).

Finally, suppose $O_{p'}(G)=1$ and $\ell_p(G)=2$. As in the proof of
Theorem \ref{55}, we see that
$\mathcal{A}_p(G)_{>X}= \mathcal{A}_p(P_1)_{>X}
\cup \ldots \cup \mathcal{A}_p(P_r)_{>X}$, with
$P_1, \ldots, P_r$ suitably chosen Sylow $p$-subgroups of $G$, and
$\mathcal{A}_p(P_i)_{>X} \cap \mathcal{A}_p(P_j)_{>X} =
\mathcal{A}_p(O_p(G))_{>X}$. Also,
$[P_i:O_p(G)]=p$. Then, by again applying section $3$ and the Gluing Lemma
as in the proof of \ref{55}, the proof is now complete. $\square$

Thus Theorem \ref{main} and Theorem \ref{example} are fully proved.


\section*{Acknowledgments}
This work is part of a Tesi di Laurea at the University of Firenze.
I would like to thank my supervisor Carlo Casolo for proposing this problem to me and for all
his help and his suggestions. I want to thank the referee for the many important comments
that have greatly improved the presentation of this paper, and for his great patience while carefully reading the work.
I would also like to thank Ken Brown
for providing helpful references and Francesco Fumagalli for proposing to me a 
prolific approach in extending my initial results.

\end{document}